\documentclass[10pt]{article}
\usepackage{amssymb,amsmath,amscd,amsthm,xcolor}
\usepackage{authblk}
\usepackage{float}
\usepackage{graphicx}
\usepackage[cm]{fullpage}

\newcommand{\Keywords}[1]{\par\noindent
{\small{\em Keywords\/}: #1}}

\usepackage{verbatim}

\let\d=\delta

\let\le=\ls
\let\ge=\gs

\def\cov{\mathop{\rm cov}\nolimits}
\newcommand{\fy}{\varphi}
\newcommand{\tend}[1]{\longrightarrow}
\newcommand{\tends}[1]{\xrightarrow[#1]{}}
\newcommand{\tendsd}{\xrightarrow{\ d\ }}
\newcommand{\abs}[1]{\lvert #1\rvert}

\newcommand{\calB}{\mathcal{B}}

\newcommand{\Nd}{\mathbb{N}}

\newcommand{\Rd}{\mathbb{R}}
\newcommand{\Zd}{\mathbb{Z}}

\newcommand{\Un}{\mathbf{1}}

\DeclareMathOperator{\Prob}{P} \DeclareMathOperator{\Mean}{E}

\renewcommand{\d}{\mathrm{d}}
\newcommand{\e}{\mathrm{e}}

\let\0=\varnothing
\let\1=\Un

\renewcommand{\le}{\leqslant}
\renewcommand{\ge}{\geqslant}

\newcommand{\edoc}{\end{document}}

\theoremstyle{plain}
\newtheorem{theorem}{\sc Theorem}[section]

\newtheorem{proposition}{\sc Proposition}[section]

\theoremstyle{definition}

\numberwithin{equation}{section}

\begin{document}
\date{}

\title{A short note on a class of statistics for estimation of the Hurst index of fractional Brownian motion}
\author[1,*,$\dag$]{K. Kubilius}
\author[1,2,**]{V. Skorniakov}

\affil[1]{Vilnius University, Institute of Mathematics and
Informatics, Akademijos 4, LT-08663, Vilnius, Lithuania}
\affil[2]{Vilnius University, Faculty of Mathematics and
Informatics, Naugarduko 24, LT-03225, Vilnius, Lithuania}

\maketitle

\let\oldthefootnote\thefootnote
\renewcommand{\thefootnote}{\fnsymbol{footnote}}
\footnotetext[1]{Corresponding author. E-mail:
\texttt{kestutis.kubilius@mii.vu.lt}} \footnotetext[2]{This research was funded by a grant (No.
MIP-048/2014) from the Research Council of Lithuania.} \footnotetext[7]{E-mail:
\texttt{viktor.skorniakov@mif.vu.lt}}
\let\thefootnote\oldthefootnote

\abstract{We propose some class of statistics suitable for
estimation of the Hurst index of the fractional Brownian motion
based on the second order increments of an observed discrete
trajectory.\Keywords{fractional Brownian motion, Hurst index,
consistent estimator, central limit theorem} }

\section{Introduction}

The aim of this short note is to present some class of statistics
suitable for estimation of the Hurst index from discretely observed
trajectory of the fractional Brownian motion (fBm). The idea behind
construction is fairly simple and makes use of self-similarity and
stationarity of the second order increments. Employment of these
properties immediately enables to prove usual asymptotic results,
namely strong consistency and normality. Perhaps the most
interesting feature of the present work is an intersection with
results of \cite{Surgailis-11} suggesting one of possible
generalizations. It appears that in case of the fBm increment ratio
(IR) statistic of \cite{Surgailis-11} belongs to the class of
statistics considered in the paper. Therefore it seems that
following along the lines of \cite{Surgailis-11} one can build a
class of statistics similar to the ones considered here and suitable
for measuring the roughness of random paths considered in
\cite{Surgailis-11}.

The structure of the paper is as follows. In section \ref{s:results}
we state theoretical result. In section \ref{s:examples} we give two
concrete examples from the class of suggested statistics. Finally
section \ref{s:proofs} contains a proof of the main theoretical
result together with a short subsection of auxiliary results needed
for the proof and collected only for the reader's convenience. The
reader unfamiliar with a topic should consult that subsection first
and then proceed to the main result.

\section{Main result}\label{s:results}

\subsection{Statement}
We assume that an observed discrete sample corresponds to the
uniform partition of a time interval of a trajectory
$(B^H_{t})_{t\in[0;T]}$ with fixed $T>0$. Since $B^H$ is
self-similar, w.l.o.g. in the rest of the paper we concentrate on
samples $B^{H}_{\frac{i}{n}}$, $i=0,\dots,n$, corresponding to a
trajectory $(B^H_{t})_{t\in[0;1]}$. Whenever it is possible, we omit
superscript and write $B_t$ instead of $B_t^H$.

Let
$d_{n,i}=d_i=\Delta^{(2)}B_{\frac{i+1}{n}}=B_{\frac{i+1}{n}}-2B_{\frac{i}{n}} +B_{\frac{i-1}{n}}$, $i=1,\dots,n-1$ be an array of the second order differences obtained from the
sample and $r_i=\frac{d_{i+1}}{d_i}$, $i=1,\dots,n-1$. Our main result
is contained in the proposition given below.

\begin{proposition}\label{p:main_result}
Let $\rho(x)=\frac{-7-9^{x}+4^{x+1}}{2(4-4^x)}$, $x\in(0;1)$, $K$
denotes a standard Cauchy r.v.\footnote{i.e. $K$ is absolutely
continuous and has a density with respect to Lebesgue measure on
$\calB(\Rd)$ given by $f_K(x)=\frac{1}{\pi(1+x^2)}$, $x\in\Rd$.} and
$h:\Rd\to\Rd$ be measurable. Assume that:
\begin{itemize}
  \item[(i)] $\Mean h^2(K+\rho(H))<\infty$;
  \item[(ii)] $H\mapsto \Mean h(K+\rho(H))$ possesses non-zero derivative of constant sign in a neighborhood of $H$.
\end{itemize}
Then
\begin{gather*}
    \widehat{H}_{n,h}=\fy\left(\bar{h}_n\right)\to H\ \text{a.s.}\quad\text{ and
    }\quad \sqrt{n}\left(\widehat{H}_{n,h}-H\right)\tendsd N(0;\sigma^2_h),\
    \text{ where}
    \quad\bar{h}_n=\frac{1}{n}\sum_{i=1}^n h(r_i),
\end{gather*}
$\fy$ denotes an inverse of $H\mapsto\Mean h(K+\rho(H))$ and
$\sigma_h^{2}$ is precisely defined in subsection \ref{ss:proof}.
\end{proposition}

\section{Concrete examples}\label{s:examples}
In this section we give two examples of functions which satisfy
conditions $(i)-(ii)$ stated in proposition \ref{p:main_result}. The
first function is considered because it usually happens that
$\arcsin$ transform symmetrizes distribution and improves normal
approximation. The second one demonstrates a connection with
\cite{Surgailis-11} discussed in the introduction. It is worthwhile
to mention that in \cite{Surgailis-11} the reader can find an
example of applications for estimation problems within a framework
of diffusions. Statistics introduced in this paper may be applied in
a similar way.

In case of the first example we check conditions $(i)-(ii)$ and give
an expression for $\fy$. In case of the second example we simply
state a form of $h$ and $\Mean h(K+\rho(H))$ referring for the
details to \cite{Surgailis-11}.

Before proceeding to the mentioned examples note that
$x\mapsto\rho(x)$ is increasing with a range equal to
$\left(-\frac{2}{3};-2+\frac{9}{8}\frac{\ln 9}{\ln4}\right)$ and
derivative
\begin{equation}\label{e:r_prime}
    \rho'(x)=\frac{1}{2}\left(\frac{-9^x\cdot4\cdot\ln(9)+4^x\cdot9\cdot\ln(4)+36^x\ln(9/4)}{(4-4^x)^2}\right).
\end{equation}
In the rest of this subsection we omit an argument for $\rho(H)$
when it appears unnecessary and write $\rho$ instead.

\medskip\emph{Example 1.} Let $h(x)=\sin x$. Then $h$ is bounded and
$(i)$ holds. Next, note that $K$ is symmetric r.v. Hence, its
characteristic function $\psi_K(t)=\Mean\cos(tK)=\e^{-\abs{t}}$ and
for any odd $g$ it holds that $\Mean g(K)=0$. Therefore
\begin{equation*}
    \Mean \sin(K+\rho)=\cos\rho\Mean\sin K+\sin\rho\Mean\cos
    K=(\sin\rho)\psi_K(1)=\frac{\sin\rho}{\e}.
\end{equation*}
Properties of $H\mapsto\rho(H)$ imply that
$H\mapsto\e^{-1}\sin\rho(H)$ is increasing on $(0;1)$ with an
inverse
\[
\fy(y)=\rho^{-1}(\arcsin(\e
y)),\qquad y\in\left(\frac{\sin(-2/3)}{\e};\frac{\sin\left(-2+\frac{9}{8}\frac{\ln
9}{\ln4}\right)}{\e}\right)
\]
and derivative
$\fy^{\prime}(\Mean\sin(K+\rho(H)))={\e}\left({\rho^{\prime}(H)\cos\rho(H)}\right)^{-1}.$

\medskip\emph{Example 2.} Setting $h(x)=\frac{\abs{1+x}}{1+\abs{x}}$
one gets statistic of \cite{Surgailis-11} with
\begin{equation*}
    \Mean
    h(K+\rho)=\frac{1}{\pi}\left(\arccos(-\rho)+\sqrt{\frac{1+\rho}{1-\rho}}\ln\left(\frac{2}{1+\rho}\right)\right).
\end{equation*}

\section{Auxiliary facts and the proof}\label{s:proofs}
\subsection{Auxiliary facts}\label{ss:auxiliary}

Below we list several facts needed for the proof of the main result.

\subsubsection{Properties of the fBm}
\begin{itemize}
  \item Fractional Brownian motion $(B_t^H)_{t\ge 0}$ is a centered continuous Gaussian process
  with a covariance function
    \begin{equation*}
    \Mean  B^H_tB^H_s=\frac{1}{2}\left(\abs{t}^{2H}+\abs{s}^{2H}-\abs{t-s}^{2H}\right),\qquad H\in(0;1).
    \end{equation*}
  \item A sequence of the second order increments  $X_i=\Delta^{(2)}B^H_{i+1}=B^H_{i+1}-2B^H_{i}+B^H_{i-1}$, $i\ge
  1$ is stationary, $\forall i\ X_{i}\sim
  N(0;4-4^{H})$,  $\mathrm{Corr}(X_{i},X_{i+1})=\rho(H)=\frac{-7-9^{H}+4^{H+1}}{2(4-4^H)}\,$.
  Moreover (see, e.g. \cite{Coeurjolly-01}),
    \begin{equation}\label{e:correl}
       \rho_k=\rho_k(H)=\mathrm{Corr}(X_{i},X_{i+k})=O\left(\frac{1}{k^{4-2H}}\right),\qquad k\to\infty.
    \end{equation}
   \item fBm is self-similar, that is $(B_{at}^H)_{t\ge 0}\stackrel{d}{=}(a^H B_t^H)_{t\ge0}$.
\end{itemize}

\subsubsection{Central limit theorem for a stationary sequence of
Gaussian random vectors} Let $(Z_i)_{i\ge
1},Z_i=(Z_{i,1},\dots,Z_{i,d})^T$ be a stationary sequence of
centered $\Rd^d$ valued Gaussian r.vs. For $k\in\Zd$ and
$p,q\in\{1,\dots,d\}$ denote
\begin{equation*}
    r^{(p,q)}(k)=\Mean Z_{m,p} Z_{m+k,q},
\end{equation*}
where $m$ is any natural number satisfying $m,m+k\ge 1$. Note that
\begin{equation}\label{e:rpq0}
    r^{(p,q)}(0)=\Mean Z_{m,p} Z_{m,q}=\Mean Z_{1,p} Z_{1,q}
\end{equation}
and
\begin{equation}\label{e:rpqk}
   r^{(p,q)}(-k){=}\Mean Z_{m,p} Z_{m-k,q}=
    \Mean Z_{k+1,p} Z_{(k+1)-k,q}=\Mean Z_{k+1,p}
    Z_{1,q}=r^{(q,p)}(k),\quad\text{ for }k\ge 1.
\end{equation}
We make use of the following result given in \cite{Arcones-94}
(Theorem 2).

\begin{theorem}\label{t:CLT_Arcones}
Assume that $f:\Rd^d\to\Rd$ is measurable, $\Mean f^2(Z_1)<\infty$
and for each $(p,q)\in\{(i,j)\mid i,j=1,\dots d\}$ there exist
finite limits
\[
    \lim_{n\to\infty}\frac{1}{n}\sum_{j,k=1}^{n}r^{(p,q)}({j-k}),\qquad
    \lim_{n\to\infty}\frac{1}{n}\sum_{j,k=1}^{n}(r^{(p,q)}({j-k}))^2.
\]
Then
\begin{equation}\label{e:normal_limit}
    \frac{1}{\sqrt{n}}\sum_{i=1}^{n}(f(Z_i)-\Mean f(Z_1))\tendsd
    N(0;\sigma^2_f),
\end{equation}
where
\[
\sigma^2_f=\mathrm{Var}(f(Z_1))+2\sum_{k=1}^{\infty}\cov(f(Z_1),f(Z_{1+k})).
\]
\end{theorem}

\subsection{Proof}\label{ss:proof}

Retain the notions introduced in the previous sections and consider
a bivariate Gaussian sequence
$Z_i=(Z_{i,1},Z_{i,2})=\frac{1}{\sqrt{4-4^H}}(X_{i+1},X_i)$, $i\ge 1$.
It follows from the properties listed above that $(Z_i)_{i\ge 1}$ is
stationary and that
\begin{equation}\label{e:bound_on_r}
    \abs{r^{(p,q)}(k)}=\abs{\Mean Z_{1,p} Z_{1+k,q}}\le \frac{C}{k^{4-2H}}
\end{equation}
for all $p,q=1,2$, $k\ge 1,$ and finite positive constant $C$ depending
only on $H$. The bound implies convergence of
$\sum_{k=1}^{\infty}r^{(p,q)}(k)$ and
$\sum_{k=1}^{\infty}(r^{(p,q)}(k))^2$. By
\eqref{e:rpq0}--\eqref{e:rpqk}
\begin{equation*}
    \frac{1}{n}\sum_{j,k=1}^n
    (r^{(p,q)}(j-k))^{i}=(r^{(p,q)}(0))^i+\frac{1}{n}\left(\sum_{j=1}^{n-1}\sum_{k=j+1}^{n}
    (r^{(q,p)}(k-j))^{i}+\sum_{j=1}^{n-1}\sum_{k=j+1}^{n}
    (r^{(p,q)}(k-j))^{i}\right),\qquad i =1,2.
\end{equation*}
Thus, it suffices to show that for any $p,q=1,2,$ there exist limits
\[
\lim_{n\to\infty}n^{-1}\sum_{j=1}^{n-1}\sum_{k=j+1}^{n}(r^{(p,q)}(k-j))^{i},\qquad i=1,2.
\]
Let $i=1$. Fix $p,q$ and note that
\begin{align*}
    \frac{1}{n}\sum_{j=1}^{n-1}\sum_{k=j+1}^{n}r^{(p,q)}({k-j})=&[k-j=l]=
    \frac{1}{n}\sum_{j=1}^{n-1}\sum_{l=1}^{n-j}r^{(p,q)}(l)=
    \sum_{l=1}^{n-1}r^{(p,q)}(l)\left(1-\frac{l}{n}\right)\\
    =&\left[\1_{\{1,\dots,n-1\}}(l)r^{(p,q)}(l)\left(1-\frac{l}{n}\right)=\psi_n(l)\right]=\sum_{l=1}^{\infty}\psi_n(l)=\int_{\Nd}\psi_n\d\mu,
\end{align*}
where $\mu$ denotes a counting measure on $\Nd$. For each fixed $l$
it holds true $\psi_n(l)\tends{n\to\infty}r^{(p,q)}(l)$. By
\eqref{e:bound_on_r},
\[\abs{\psi_n(l)}\le\psi(l)=\frac{C}{l^{4-2H}}\,,\qquad l\in\Nd.\] Since $\psi$
is integrable with respect to $\mu$, Dominated Convergence theorem yields relationship
\begin{equation*}
    \lim_{n\to\infty}\int_{\Nd}\psi_n\d\mu=\int_{\Nd}\lim_{n\to\infty}\psi_n\d\mu=\sum_{l=1}^{\infty}r^{(p,q)}(l)\in\Rd.
\end{equation*}
Hence,
\[
\lim_{n\to\infty}n^{-1}\sum_{j=1}^{n-1}\sum_{k=j+1}^{n}r^{(p,q)}(k-j)=\sum_{l=1}^{\infty}r^{(p,q)}(l).
\]
Identical argument shows that
\begin{equation*}
    \lim_{n\to\infty}n^{-1}\sum_{j=1}^{n-1}\sum_{k=j+1}^{n}\big(r^{(p,q)}(k-j)\big)^{2} =\sum_{l=1}^{\infty}\big(r^{(p,q)}(l)\big)^2.
\end{equation*}
Therefore theorem \ref{t:CLT_Arcones} applies to $(Z_i)$ provided
function $f$ is suitably chosen. Let
\begin{equation*}
    f(x,y)=\1_{\Rd\times(\Rd\setminus\{0\})}(x,y)h\left(\frac{x}{y}\right),\qquad
    R_i=\frac{Z_{i,1}}{Z_{i,2}}\,, \quad i\ge1.
\end{equation*}
Recall that a ratio of two independent standard Gaussian r.vs. has
the standard Cauchy distribution. Keeping this in a view we arrive
to the following conclusions\footnote{$K$ denotes a r.v. having
standard Cauchy distribution (see subsection \ref{s:results})} :
\begin{itemize}
  \item $R_i=\frac{Z_{i,1}}{Z_{i,2}}=\frac{(X_{i+1}-\rho(H)X_{i})/\sqrt{4-4^H}}{X_{i}/\sqrt{4-4^H}}+\rho(H)\stackrel{d}{=}
  K+\rho$, since $(X_{i+1}-\rho(H)X_{i})/\sqrt{4-4^H}$, ${X_{i}/\sqrt{4-4^H}}\sim N(0;1)$ are uncorrelated;
  \item by self-similarity $(R_i)\stackrel{d}{=}(r_i)$;
  \item because of absolute continuity of
Gaussian r.v. and the previous facts $\Mean f^j(Z_{1})=\Mean
h^j(R_1)=\Mean h^j(r_1)=\Mean h^j(K+\rho(H)),$ $j=1,2$.
\end{itemize}
Also note that in addition to stationarity
$(X_i)_{i\ge 1}$ has vanishing correlation (see (\ref{e:correl})).
Therefore it is ergodic and conclusions listed above together with
Ergodic theorem (see Corollary 8.6.3 in \cite{Athreya-2006}) yield
relationships
\begin{equation*}
    1=\Prob\left(\lim_{n\to\infty}\frac{f(Z_1)+\dots+f(Z_n)}{n}=\Mean
    f(Z_1)\right)=\Prob\left(\lim_{n\to\infty}\frac{h(r_1)+\dots+h(r_n)}{n}=\Mean
    h(r_1)=\Mean h(R_1)\right).
\end{equation*}
That is, $\bar h_n\to \Mean h(R_1)$ a.s. and by continuous mapping
theorem $\fy(\bar h_n)\to\fy( \Mean h(R_1))=\fy(\Mean
h(K+\rho(H)))=H$ a.s.

Next, note that from previously stated equality $\Mean
f^2(Z_{1})=\Mean h^2(K+\rho(H))$ and condition $\Mean
h^2(K+\rho(H))<\infty$ it follows that Theorem \ref{t:CLT_Arcones}
applies to the chosen $f$ and gives \eqref{e:normal_limit} which may
be rewritten as
\begin{gather}\label{e:normal_limit2}
    \frac{1}{\sqrt{n}}\sum_{i=1}^{n}\big(h(r_i)-\Mean h(r_1)\big)=
    {\sqrt{n}}\big(\bar{h}_n-\Mean h(r_1)\big)\tendsd
    N(0;\sigma^2_f),\\
    \sigma^2_f=\mathrm{Var}(h(R_1))+2\sum_{k=1}^{\infty}\cov(h(R_1),h(R_{1+k})).\nonumber
\end{gather}
To finish the proof, one has to apply the Delta method, which also
yields asymptotic variance $\sigma_h^2=(\fy^\prime(\Mean
h(K+\rho(H))))^2\sigma^2_f$. \qed


\begin{thebibliography}{00}

\bibitem{Arcones-94} M. A. Arcones,
Limit theorems for nonlinear functionals of a stationary Gaussian
sequence of vectors, \emph{The Annals of Probability}, \textbf{22}(4) (1994), 2242-–2274.

\bibitem{Athreya-2006} Athreya, Krishna B. and Lahiri, Soumen N.,
Measure Theory and Probability Theory (Springer Texts in Statistics), 2006, Springer-Verlag New York, Inc.


\bibitem{Coeurjolly-01}  J.-F. Coeurjolly, Estimating the parameters of a fractional Brownian motion by discrete variations of its sample paths,
\emph{Statistical Inference for Stochastic Processes}, \textbf{4} (2001), 199--227.


\bibitem{Surgailis-11}  J. M. Bardet and D. Surgailis, Measuring the roughness of random paths by increment ratios, \emph{Bernoulli}, \textbf{17}(2) (2011), 749--780.


\end{thebibliography}
\end{document}